\documentclass{amsart}
\usepackage{amsfonts, amssymb, latexsym, amsmath, amscd, 
graphics}
\usepackage[mathscr]{eucal}
\usepackage{enumerate}

\theoremstyle{plain}
\newtheorem{theorem}{Theorem}[section]

\newtheorem{lemma}[theorem]{Lemma}
\newtheorem{proposition}[theorem]{Proposition}
\newtheorem{sublemma}[theorem]{Sublemma}

\theoremstyle{definition}
\newtheorem{definition}[theorem]{Definition}

\newtheorem{remark}[theorem]{Remark}

\newcommand{\mbf}[1]{\boldsymbol{\mathit{#1}}}
\newcommand{\wh}[1]{\widehat{#1}}
\newcommand{\wt}[1]{\widetilde{#1}}
\newcommand{\mscr}[1]{\mathscr{#1}}
\newcommand{\mcal}[1]{\mathcal{#1}}
\newcommand{\pref}[1]{(\ref{#1})}

\newcommand{\bpref}[1]{\textbf{(\ref{#1})}}

\newcommand{\parop}[1]{\partial_{#1}}
\newcommand{\parfracop}[1]{\frac{\partial}{\partial #1}}

\newcommand{\psup}[2]{#1^{(#2)}}
\newcommand{\parfrac}[2]{\frac{\partial #1}{\partial #2}}
\newcommand{\subsup}[3]{#1_{#2}^{#3}}

\newcommand{\subpsup}[3]{#1_{#2}^{(#3)}}

\newcommand{\NN}{\mathbb{N}}
\newcommand{\ZZ}{\mathbb{Z}}
\newcommand{\RR}{\mathbb{R}}
\newcommand{\leqs}{\leqslant}
\newcommand{\geqs}{\geqslant}

\newcommand{\ve}{\varepsilon}

\newcommand{\vp}{\varphi}

\newcommand{\al}{\alpha}
\newcommand{\be}{\beta}
\newcommand{\ga}{\gamma}
\newcommand{\de}{\delta}

\newcommand{\ze}{\zeta}
\newcommand{\thet}{\theta}

\newcommand{\ka}{\kappa}
\newcommand{\la}{\lambda}
\newcommand{\rh}{\rho}
\newcommand{\si}{\sigma}
\newcommand{\ta}{\tau}

\newcommand{\om}{\omega}

\newcommand{\Ga}{\Gamma}
\newcommand{\De}{\Delta}

\newcommand{\La}{\Lambda}

\newcommand{\Om}{\Omega}

\newcounter{hrealm}
\newcounter{prealm}
\newcounter{grealm}
\newcounter{sarealm}

\renewcommand{\thehrealm}{H\arabic{hrealm}}
\renewcommand{\theprealm}{P\arabic{prealm}}
\renewcommand{\thegrealm}{G\arabic{grealm}}
\renewcommand{\thesarealm}{SA\arabic{sarealm}}

\numberwithin{equation}{section}

\begin{document}

\title[Strange attractors and degenerate Hopf bifurcations]{Strange attractors
in periodically-kicked\\
degenerate Hopf bifurcations}

\author{William Ott}
\address{Courant Institute of Mathematical Sciences\\
New York, New York 10012}
\email[William Ott]{ott@cims.nyu.edu}
\urladdr[William Ott]{www.cims.nyu.edu/$\sim$ott}

\keywords{Degenerate Hopf bifurcation, Rank one map, Misiurewicz map}

\subjclass[2000]{Primary:}

\date{June 2007}

\begin{abstract}
We prove that spiral sinks (stable foci of vector fields) can be transformed
into strange attractors exhibiting sustained, observable chaos if subjected to
periodic pulsatile forcing.  We show that this phenomenon occurs in the context
of periodically-kicked degenerate supercritical Hopf bifurcations.  The results
and their proofs make use of a new $k$-parameter version of the theory of rank
one maps developed by Wang and Young.
\end{abstract}

\maketitle

\section{Introduction}\label{s:intro}

This paper aims to support the idea that shear and twist are natural mechanisms
for the production of sustained, observable chaos in forced dynamical systems.
Consider a weakly stable dynamical structure such as an equilibrium point or a
limit cycle.  If shear or twist is present, then forcing of various types can
transform the weakly stable structure into a strange attractor.  The nature of
the forcing is not essential.  Admissible types of forcing include periodic
pulsatile drives, deterministic continuous-time signals, and random signals
generated by stochastic processes.  The strange attractors possess many of the
dynamical, statistical, and geometrical properties commonly associated with
chaotic dynamics.

We study the simplest weakly stable dynamical structure.  This is the spiral
sink (or stable focus), an equilibrium point of a vector field with the property
that the linearization of the field at the equilibrium point has a pair of
complex conjugate eigenvalues $\al \pm i \be$ satisfying $\al < 0$ and $\be \neq
0$.  We consider the degenerate supercritical Hopf bifurcation in two
dimensions.  When a generic supercritical Hopf bifurcation occurs, the spiral
sink becomes unstable and a limit cycle is born.  In the degenerate case, the
spiral sink loses its stability but no limit cycle is born.  We prove that in
the case of the degenerate supercritical Hopf bifurcation, periodic pulsatile
drives transform the spiral sink into a strange attractor.  The analysis is
certainly not limited to Hopf bifurcations.  We work in this context because the
origin of the shear is transparent in the defining differential equations.

The analysis is based on the beautiful dynamical theory of rank one maps
formulated by Wang and Young~\cite{WqYls2001, WqYls2005}.  Speaking
impressionistically, rank one maps are strongly dissipative maps exhibiting a
single direction of instability.  Rank one theory provides checkable conditions
that imply the existence of strange attractors for a positive-measure set of
parameters within a given parametrized family of rank one maps.  The conditions
appear within the following scheme.
\begin{enumerate}[(1)]
\item
\label{li:pro1}
Let dissipation go to infinity.  This procedure produces the singular limit, a
parametrized family of one-dimensional maps.
\item
\label{li:pro2}
Check that the singular limit includes a map with strong expanding properties (a
map of Misiurewicz type).
\item
\label{li:pro3}
Verify a parameter transversality condition.
\item
\label{li:pro4}
Verify a nondegeneracy condition.  This allows information about the singular
limit to be passed to the maps with finite dissipation.
\end{enumerate}
Steps~\pref{li:pro3} and~\pref{li:pro4} cumulatively require verifying that only
finitely many quantities do not vanish.  For good parameters, parameters
corresponding to maps admitting strange attractors, rank one theory provides a
reasonably complete dynamical description of the map.

The attractor supports a positive, finite number of ergodic SRB measures.  The
orbit of Lebesgue almost-every point in the basin of attraction has a positive
Lyapunov exponent and is asymptotically distributed according to one of the
ergodic SRB measures.  Each SRB measure satisfies the central limit theorem and
exhibits exponential decay of correlations.  A symbolic coding exists for orbits
on the attractor.  This symbolic coding implies the existence of equilibrium
states and a measure of maximal entropy.  Summarizing, the map has a
nonuniformly hyperbolic character and exhibits sustained, observable chaos.

Of the four steps in the rank one scheme, step~\pref{li:pro2} is the most
fundamental and typically requires the most work.  A Misiurewicz map has the
property that the positive orbit of {\bfseries \itshape every} critical point
remains bounded away from the critical set.  Existing papers on rank one theory
view the singular limit as a one-parameter family $\{ f_{a} \}$ of
one-dimensional maps.  This view makes locating Misiurewicz parameters difficult
if the maps have multiple critical points.  If each map $f_{a}$ has exactly one
critical point $c(a)$, then locating Misiurevicz parameters is relatively easy.
Assuming that $f_{a} (c(a))$ moves reasonably quickly as one varies $a$, simply
locate an invariant set $\La (a)$ that is disjoint from the critical set and
then choose $a^{*}$ such that $f_{a^{*}} (c(a^{*})) \in \La (a^{*})$.  The set
$\La (a)$ could be a periodic orbit or a Cantor set.  If the singular limit
consists of maps with multiple critical points, then one must locate a parameter
$a^{*}$ for which all of the critical orbits of $f_{a^{*}}$ are contained in
good invariant sets.  This is a serious challenge because good invariant sets
such as periodic orbits and Cantor sets typically have Lebesgue measure zero.
Wang and Young~\cite{WqYls2002} overcome this challenge.  However, their results
assume that the maps in the singular limit possess an extremely large amount of
expansion.

We prove that significantly less expansion is needed if the singular limit is
viewed as an $m$-parameter family for $m$ sufficiently large.  Assume that the
singular limit consists of maps with $k$ critical points.  We prove that if the
singular limit is viewed as a $k$-parameter family, then it contains Misiurewicz
points assuming the maps are mildly expanding and assuming the parameters are
independent in a sense to be made precise.  This result widens the scope of rank
one theory.

We view this work as an element of a growing list of applications of rank one
theory.  The theory has been rigorously applied to simple mechanical
systems~\cite{WqYls2002}, periodically-kicked limit cycles and Hopf
bifurcations~\cite{WqYls2003}, and the Chua circuit~\cite{WqOa2005}.
Guckenheimer, Wechselberger, and Young~\cite{GjWmYls2006} connect rank one
theory and geometric singular perturbation theory by formulating a general
technique for proving the existence of chaotic attractors for three-dimensional
vector fields with two time scales.  Lin~\cite{Lk2006} demonstrates how rank one
theory can be combined with sophisticated computational techniques to analyze
the response of concrete nonlinear oscillators of interest in biological
applications to periodic pulsatile drives.  Lin and Young~\cite{LkYls2007} study
shear-induced chaos numerically in situations beyond the reach of current
analytical tools.  In particular, they consider stochastic forcing.  This work
supports the belief that shear-induced chaos is both widespread and robust.

We organize the presentation of ideas as follows.  In Section~\ref{s:pkdhb}, we
present the main results for periodically-kicked degenerate Hopf bifurcations.
In Section~\ref{s:mis}, we prove the result concerning the existence of
Misiurewicz points in $k$-parameter families of one-dimensional maps and we
present a two-parameter example.  Section~\ref{s:r1t} presents rank one theory
viewing the singular limit as a $k$-parameter family.  Finally, in
Section~\ref{s:pst} we prove the results presented in Section~\ref{s:pkdhb}.

\section{Periodically-kicked degenerate Hopf bifurcations}\label{s:pkdhb}

The normal form for the supercritical Hopf bifurcation in two spatial dimensions
is given in polar coordinates by
\begin{equation*}
\left\{
\begin{aligned}
\dot{r} &= (\mu - \al_{\mu} r^{2})r + r^{5} g_{\mu} (r, \thet)\\
\dot{\thet} &= \om + \ga_{\mu} \mu + \be_{\mu} r^{2} + r^{4} h_{\mu} (r, \thet)
\end{aligned}
\right.
\end{equation*}
Here $\mu$ is the bifurcation parameter and $\om$ is a constant.  The
multipliers $\al_{\mu}$, $\ga_{\mu}$, and $\be_{\mu}$ depend smoothly on $\mu$.
The functions $g_{\mu}$ and $h_{\mu}$ depend smoothly on $\mu$ and they are of
class $C^{4}$ with respect to $r$ and $\thet$.  The normal form for the
degenerate Hopf bifurcation in two spatial dimensions is obtained by setting
$\al_{\mu} = 0$ for all $\mu$ and replacing $\mu$ with $-\mu$, yielding
\begin{equation}\label{e:dhbnf}
\left\{
\begin{aligned}
\dot{r} &= - \mu r + r^{5} g_{\mu} (r, \thet)\\
\dot{\thet} &= \om + \ga_{\mu} \mu + \be_{\mu} r^{2} + r^{4} h_{\mu} (r, \thet)
\end{aligned}
\right.
\end{equation}
For $\mu > 0$, the origin is an asymptotically stable equilibrium point (a
sink).  We study~\pref{e:dhbnf} in this $\mu$-range.  Let $\wh{F}_{t}$ denote
the flow generated by~\pref{e:dhbnf}.  We perturb the flow $\wh{F}_{t}$ with a
\textquoteleft kick' map $\ka$ defined as follows.  Let $L > 0$ and let $\rh_{2}
> 0$.  The map $\ka = \ka_{\mu, L, \rh_{2}}$ is given in rectangular coordinates
by
\begin{equation*}
\ka
\begin{pmatrix}
r \cos (\thet)\\
r \sin (\thet)
\end{pmatrix}
=
\begin{pmatrix}
r \cos (\thet)\\
r \sin (\thet) + L \mu^{\rh_{2}}
\end{pmatrix}
\end{equation*}
The composition $\wh{F}_{t} \circ \ka$ may be thought of as a perturbation
followed by a period of relaxation.  We define an annulus map associated with
$\wh{F}_{t} \circ \ka$.  Let $\mscr{A}$ denote the annulus defined by
\begin{equation*}
\mscr{A} = \{ (r, \thet) : K_{4}^{-1} \mu^{\rh_{1}} \leqs r \leqs K_{4}
\mu^{\rh_{1}} \},
\end{equation*}
where $K_{4} > 1$ and $0 < \rh_{2} < \rh_{1}$.  Let $\tilde{r}$ denote the
distance from $\ka (\mscr{A})$ to the origin.  We have $\tilde{r} = L
\mu^{\rh_{2}} - K_{4} \mu^{\rh_{1}}$.  Define the relaxation time $\ta (\mu)$ by
\begin{equation*}
\tilde{r} e^{-\mu \ta (\mu)} = \mu^{\rh_{1}}.
\end{equation*}
For $\mu$ sufficiently large, $\wh{F}_{\ta (\mu)} \circ \ka$ maps $\mscr{A}$ into
$\mscr{A}$.  

The following theorem states that under certain conditions, the annulus map
$\wh{F}_{\ta (\mu)} \circ \ka$ admits a strange attractor for a positive-measure
set of values of $\mu$.  We make the crucial assumption that the twist factor
$\be_{0}$ is nonzero.  A nonzero twist factor implies the existence of an
angular-velocity gradient in the radial direction for values of $\mu$ in a
neighborhood of the bifurcation parameter $\mu = 0$.  This angular-velocity
gradient allows the flow to stretch and fold the phase space, thereby producing
chaos.

The chaos in this setting is sustained in time and observable.  The strange
attractors possess many of the geometric and dynamical properties normally
associated with chaotic systems.  These properties include the existence of a
positive Lyapunov exponent~\bpref{li:sa1}, the existence of SRB measures and
basin property~\bpref{li:sa2}, and statistical properties such as exponential
decay of correlations and the central limit theorem for dynamical
observations~\bpref{li:sa3}.  In addition, if $\frac{L}{|\be_{0}|}$ is
sufficiently large, then the annulus map admits a unique SRB
measure~\bpref{li:sa4}.  Properties~\bpref{li:sa1}-\bpref{li:sa4} are described
in detail in Section~\ref{s:r1t}.

\begin{theorem}\label{t:sink}
Assume $\be_{0} \neq 0$.  Let $\rh_{1}$ and $\rh_{2}$ satisfy $\rh_{2} \in
(\frac{3}{8}, \frac{1}{2})$ and $\rh_{1} + \rh_{2} = 1$.
\begin{enumerate}[(1)]
\item
\label{li:sink1}
There exists $M_{0} > 0$ such that if $L \geqs \frac{M_{0}}{|\be_{0}|}$, then
there exist $L^{*} \in [L, L + \frac{\pi}{|\be_{0}|}]$ and $\mu_{0} > 0$
satisfying the following.  The parameter interval $(0, \mu_{0}]$ contains a set
$\De = \De (L^{*})$ of positive measure such that for $\mu \in \De$, the map
$\wh{F}_{\ta (\mu)} \circ \ka$ admits a strange attractor with
properties~\bpref{li:sa1},~\bpref{li:sa2}, and~\bpref{li:sa3}.  The set $\De$
intersects every interval of the form $(0, \tilde{\mu}]$ in a set of positive
measure.
\item
\label{li:sink2}
There exists $M_{1} \gg M_{0}$ such that for all $L \geqs
\frac{M_{1}}{|\be_{0}|}$, there exists a set $\De = \De (L)$ with the properties
described in~\pref{li:sink1}.
\item
\label{li:sink3}
If $L$ is sufficiently large and $\mu \in \De (L)$, then~\bpref{li:sa4} holds as
well.
\end{enumerate}
\end{theorem}

\section{Locating Misiurewicz points}\label{s:mis}

Let $I$ denote an interval or the circle $S^{1}$.  Let $F : I
\times [a_{1}, a_{2}] \times [b_{1}, b_{2}] \to I$ be a $C^{2}$ map.  The map
$F$ defines a two-parameter family $\mscr{F} = \{ f_{a,b} : a \in [a_{1},
a_{2}], \: b \in [b_{1}, b_{2}] \}$ via $f_{a,b} (x) = F(x,a,b)$.  Set $A =
[a_{1}, a_{2}]$ and $B = [b_{1}, b_{2}]$.  We assume that for each $(a,b) \in A
\times B$, $f_{a,b}$ has two critical points.  We label these critical points
$\psup{c}{1} (a,b)$ and $\psup{c}{2} (a,b)$.  Let $C = C(a,b) = \{ \psup{c}{1}
(a,b), \psup{c}{2} (a,b) \}$.  For $\de > 0$, let $C_{\de}$ denote the
$\de$-neighborhood of $C$ in $I$.

We seek to identify conditions under which $\mscr{F}$ contains strongly
expanding (Misiurewicz) maps.  We now introduce this class.

\begin{definition}\label{d:mis}
We say that $f \in C^{2} (I,I)$ is a Misiurewicz map and we write $f \in
\mscr{M}$ if the following hold for some neighborhood $V$ of $C$.
\begin{enumerate}[(A)]
\item
\textbf{(Outside of} {\mathversion{bold} $V$}\textbf{)}  There
exist $\lambda_{0} > 0$, $M_{0} \in \ZZ^{+}$, and $0 < d_{0} \leqs 1$ such
that\label{li:ma}
\begin{enumerate}[(1)]
\item
for all $n \geqs M_{0}$, if $f^{k} (x) \notin V$ for $0 \leqs k \leqs n-1$, then
$|(f^{n})' (x)| \geqs e^{\lambda_{0} n}$,\label{li:ma1}
\item
for any $n \in \ZZ^{+}$, if $f^{k} (x) \notin V$ for $0 \leqs k \leqs n-1$ and
$f^{n} (x) \in V$, then $|(f^{n})' (x)| \geqs d_{0} e^{\lambda_{0}
n}$.\label{li:ma2}
\end{enumerate}
\item
\textbf{(Critical orbits)}  For all $c \in C$ and $n > 0$, $f^{n} (c) \notin
V$.\label{li:mb}
\item
\textbf{(Inside} {\mathversion{bold} $V$}\textbf{)}\label{li:mc}
\begin{enumerate}[(1)]
\item
We have $f''(x) \neq 0$ for all $x \in V$, and\label{li:mc1}
\item
for all $x \in V \setminus C$, there exists $p_{0} (x) > 0$ such that $f^{j} (x)
\notin V$ for all $j < p_{0} (x)$ and $|(f^{p_{0} (x)})' (x)| \geqs d_{0}^{-1}
e^{\tfrac{1}{3} \lambda_{0} p_{0} (x)}$.\label{li:mc2}
\end{enumerate}
\end{enumerate}
\end{definition}

% Use of the updated user-defined commands (Greek, subsup) begins

We first formulate hypotheses that imply the existence of maps in $\mscr{F}$
that satisfy Definition~\ref{d:mis}\pref{li:mb}.

\subsection{The general result}\label{ss:gr}

We formulate the result for two-parameter families consisting of maps with two
critical points.  The result generalizes in a natural way for $k$-parameter
families consisting of maps with $k$ critical points.

The first hypothesis is formulated in terms of the evolutions
\begin{equation*}
(a,b) \mapsto \subpsup{\ga}{n}{i} (a,b) \text{ where } \subpsup{\ga}{n}{i} (a,b)
= f_{a,b}^{n} (\psup{c}{i} (a,b)).
\end{equation*}
The evolutions $\{ \subpsup{\ga}{n}{i} : n \in \NN \}$ generate {\bfseries
\itshape critical curve} dynamics.  Define $\Ga_{n} : A \times B \to
I \times I$ by $\Ga_{n} = (\subpsup{\ga}{n}{1}, \subpsup{\ga}{n}{2})$.

We now present the general hypotheses.  For $J \subset I$ and $\ve > 0$, let
$J^{\ve}$ denote the $\ve$-neighborhood of $J$.  Suppose there exist
subintervals $I_{1}$ and $I_{2}$ of $I$, subintervals $\De_{1} \subset A$ and
$\De_{2} \subset B$, $\de_{1} > 0$, and $\ve_{1} > 0$ such that the following
hold.

\begin{list}{\bfseries (\thehrealm)}
{
\usecounter{hrealm}
\setlength{\topsep}{1.5ex plus 0.2ex minus 0.2ex}
\setlength{\labelwidth}{1.2cm}
\setlength{\leftmargin}{1.5cm}
\setlength{\labelsep}{0.3cm}
\setlength{\rightmargin}{0.5cm}
\setlength{\parsep}{0.5ex plus 0.2ex minus 0.1ex}
\setlength{\itemsep}{0ex plus 0.2ex}
}
\item
\textbf{(Finite Misiurevicz condition)}  There exists $n_{0} \in \ZZ^{+}$ such
that $\Ga_{n_{0}} (\De_{1} \times \De_{2}) \supset I_{1} \times I_{2}$ and for
$i \in \{ 1, 2 \}$, $(a,b) \in \De_{1} \times \De_{2}$, and $n < n_{0}$, we have
$\subpsup{\ga}{n}{i} (a,b) \in I \setminus C_{\de_{1}} (a,b)$.
\label{li:h1}
\item
There exist fixed parameters $\hat{a} \in \De_{1}$ and $\hat{b} \in \De_{2}$
satisfying $f_{\hat{a}, \hat{b}} (I_{1}) \times f_{\hat{a}, \hat{b}} (I_{2})
\supset \subsup{I}{1}{\ve_{1}} \times \subsup{I}{2}{\ve_{1}}$.
\label{li:h2}
\item
For all $(a,b) \in \De_{1} \times \De_{2}$, we have $I_{1} \times I_{2} \subset
I \setminus C_{\de_{1}} (a,b) \times I \setminus C_{\de_{1}} (a,b)$.
\label{li:h3}
\end{list}

\begin{proposition}\label{p:m}
Suppose $\mscr{F}$ satisfies~\pref{li:h1}-\pref{li:h3}.  If
\begin{equation}\label{e:parset}
2 \sqrt{2} \max \{ \| \parop{a} F \|_{C^{0}}, \| \parop{b} F \|_{C^{0}} \} \cdot
\max \{ |\De_{1}|, |\De_{2}| \} < \ve_{1},
\end{equation}
then there exists $(a^{*}, b^{*}) \in \De_{1} \times \De_{2}$ such that for $i
\in \{ 1, 2 \}$ and for every $n \in \NN$, $\subpsup{\ga}{n}{i} (a^{*}, b^{*})
\in I \setminus C_{\de_{1}} (a^{*}, b^{*})$.
\end{proposition}

\begin{proof}[Proof of Proposition~\ref{p:m}]
Define $G = (f_{\hat{a}, \hat{b}}, f_{\hat{a}, \hat{b}})$.
Applying~\pref{li:h1} and~\pref{li:h2}, we have $G(\Ga_{n_{0}} (\De_{1} \times
\De_{2})) \supset \subsup{I}{1}{\ve_{1}} \times \subsup{I}{2}{\ve_{1}}$.  For
every $(a,b) \in \De_{1} \times \De_{2}$, we have
\begin{equation*}
\| \Ga_{n_{0} + 1} (a,b) - G(\Ga_{n_{0}} (a,b)) \| < \ve_{1}
\end{equation*}
since $\ve_{1}$ satisfies~\pref{e:parset}.  Therefore, $\Ga_{n_{0} + 1} (\De_{1}
\times \De_{2}) \supset I_{1} \times I_{2}$.  Inductively, $\Ga_{n} (\De_{1}
\times \De_{2}) \supset I_{1} \times I_{2}$ for all $n \geqs n_{0}$.  Define
\begin{equation*}
\Psi_{n_{0}} = (\De_{1} \times \De_{2}) \cap \subsup{\Ga}{n_{0}}{-1} (I_{1}
\times I_{2}).
\end{equation*}
For $n > n_{0}$, define
\begin{equation*}
\Psi_{n+1} = \Psi_{n} \cap \subsup{\Ga}{n+1}{-1} (I_{1} \times I_{2}).
\end{equation*}
Let
\begin{equation*}
\Psi = \bigcap_{k = n_{0}}^{\infty} \Psi_{k}
\end{equation*}
and choose $(a^{*}, b^{*}) \in \Psi$.
\end{proof}

\subsection{Verifying~\pref{li:h1}}\label{ss:vh1}

We present a two-step procedure for the verification of hypothesis~\ref{li:h1}.
First, we assume that $\Ga_{1}$ is a diffeomorphism on $\De_{1} \times
\De_{2}$.  This implies that the image of $\De_{1} \times \De_{2}$ contains a
rectangle in $I \times I$.  Second, if we assume that each map $f_{a,b}$ is
expanding on $I \setminus C_{\de_{1}}$, then the evolutions $\psup{\ga}{1}$ and
$\psup{\ga}{2}$ will enlarge the rectangle to macroscopic size.  The required
time for this enlargement depends upon the magnitude of the expansion.
Therefore, greater expansion results in a smaller value of $n_{0}$.  We now make
these ideas precise.

Suppose that $\Ga_{1}$ is a diffeomorphism on $\De_{1} \times \De_{2}$ such that
for $i \in \{ 1, 2 \}$ and for every $(a,b) \in \De_{1} \times \De_{2}$, we have
$\subpsup{\ga}{1}{i} (a,b) \in I \setminus C_{\de_{1}}$.  Define
\begin{equation*}
J(a,b) = \left|
\begin{matrix}
\parop{a} \subpsup{\ga}{1}{1} (a,b) & \parop{b} \subpsup{\ga}{1}{1} (a,b)\\
\parop{a} \subpsup{\ga}{1}{2} (a,b) & \parop{b} \subpsup{\ga}{1}{2} (a,b)
\end{matrix}
\right|.
\end{equation*}
Assume that there exists $k_{0} > 0$ such that $|J| \geqs k_{0}$ on $\De_{1}
\times \De_{2}$.  This implies that $\Ga_{1} (\De_{1} \times \De_{2})$ contains
a box with side length bounded below by
\begin{equation*}
\frac{\subsup{k}{0}{2}}{\la_{M}} \min \{ |\De_{1}|, |\De_{2}| \},
\end{equation*}
where
\begin{equation*}
\la_{M} = \sup_{(a,b) \in \De_{1} \times \De_{2}} \sup \{ |\la| : \la \text{ is
an eigenvalue of } D\Ga_{1}^{*} D\Ga_{1} \}.
\end{equation*}
Now suppose that for every $(a,b) \in \De_{1} \times \De_{2}$ we have
$|f_{a,b}'| \geqs K > 1$ on $I \setminus C_{\de_{1}}$.  Lower bounds on $K$ will
be given as the discussion proceeds.

We choose $K$ based on the magnitudes of the partial derivatives of
$\subpsup{\ga}{1}{1}$ and $\subpsup{\ga}{1}{2}$.  Assume there exists $\rh > 0$
such that on $\De_{1} \times \De_{2}$ we have
\begin{enumerate}
\item
$|\parop{a} \subpsup{\ga}{1}{1}| \geqs \rh$ or $|\parop{b} \subpsup{\ga}{1}{1}|
\geqs \rh$, and
\label{li:inde1}
\item
$|\parop{a} \subpsup{\ga}{1}{2}| \geqs \rh$ or $|\parop{b} \subpsup{\ga}{1}{2}|
\geqs \rh$.
\label{li:inde2}
\end{enumerate}
Suppose for the sake of definiteness that~\pref{li:inde1} and~\pref{li:inde2}
hold with respect to the operator $\parop{a}$.  We now relate spatial and
parametric derivatives.  The equation
\begin{equation}\label{e:cr}
\parfracop{a} \subpsup{\ga}{n+1}{i} (a,b) = f_{a,b}' (\subpsup{\ga}{n}{i} (a,b))
\cdot \parfracop{a} \subpsup{\ga}{n}{i} (a,b) + \parfrac{F}{a}
(\subpsup{\ga}{n}{i} (a,b), a, b)
\end{equation}
implies that parametric derivatives grow exponentially provided that spatial
derivatives grow exponentially.  If $K$ satisfies
\begin{gather*}
K \rh - \| \parop{a} F \|_{C^{0}} \geqs \frac{3}{4} K, \text{ and}\\
\| \parop{a} F \|_{C^{0}} \sum_{j=2}^{\infty} K^{-j} \leqs \frac{1}{4},
\end{gather*}
then~\pref{e:cr} implies that
\begin{equation}\label{e:dccexpg}
|\parop{a} \subpsup{\ga}{n}{i} (a,b)| \geqs \frac{1}{2} K^{n}
\end{equation}
provided $\subpsup{\ga}{j}{i} (a,b) \in I \setminus C_{\de_{1}}$ for $j < n$.
Hypothesis~\ref{li:h1} may be verified as follows.  Look for a time $n_{0}$ such
that for $i \in \{ 1, 2 \}$, $\subpsup{\ga}{n_{0}}{i} (\De_{1} \times \{ b \})
\supset I_{i}$ for every $b \in \De_{2}$ and $\subpsup{\ga}{k}{i} (a,b) \in I
\setminus C_{\de_{1}}$ for all $k < n_{0}$ and $(a,b) \in \De_{1} \times
\De_{2}$.  By~\pref{e:dccexpg}, we have
\begin{equation*}
K^{n_{0}} \approx \frac{2 \la_{M} \max \{ |I_{1}|, |I_{2}| \}}{\subsup{k}{0}{2}
\min \{ |\De_{1}|, |\De_{2}| \}}.
\end{equation*}

\subsection{A two-parameter example}\label{ss:twopex}

Let $S^{1} = \RR / 2 \pi \ZZ$.  Let $\Phi : S^{1} \to \RR$ be a $C^{3}$ function
with two nondegenerate critical points $\psup{c}{1}$ and $\psup{c}{2}$.  We
assume that $\Phi (\psup{c}{1}) \neq \Phi (\psup{c}{2})$.  Fix $\ze \in S^{1}$.
Consider the two-parameter family of circle maps $\mscr{F} = \{ f_{a,L} : a \in
S^{1}$, $L \in \RR^{+} \}$ defined by
\begin{equation*}
f_{a,L} (\thet) = \ze + L \Phi (\thet) + a.
\end{equation*}
Small perturbations of this family frequently arise as singular limits of rank
one families.  

\begin{definition}\label{d:mispair}
We say that $(a,L)$ is a {\bfseries \itshape Misiurewicz pair} if $f_{a,L} \in
\mscr{M}$.
\end{definition}

\noindent
The goal of this subsection is to prove the following result.

\begin{theorem}\label{t:al}
There exists $L_{0} > 0$ such that if $L \geqs L_{0}$, then there exists a
Misiurewicz pair $(a^{*}, L^{*})$ with $L^{*} \in [L, L + 2 \pi / |\Phi
(\psup{c}{2}) - \Phi (\psup{c}{1})|]$ and $a^{*} \in [0, 2 \pi)$.
\end{theorem}

\begin{remark}
Misiurewicz points occur with greater frequency as $L$ increases.  Wang and
Young~\cite{WqYls2002} prove that there exists $L_{1} \gg L_{0}$ such that if $L
\geqs L_{1}$, then $f_{a,L} \in \mscr{M}$ for a $\mcal{O} (1/L)$-dense subset of
parameters $a \in [0, 2 \pi)$.
\end{remark}

\begin{proof}[Proof of Theorem~\ref{t:al}]
We prove Theorem~\ref{t:al} in two steps.  We first show that for $L$
sufficiently large, if $f_{a,L}$ satisfies Definition~\ref{d:mis}\pref{li:mb},
then $f_{a,L} \in \mscr{M}$.  We then prove the existence of parameters for
which $f_{a,L}$ satisfies Definition~\ref{d:mis}\pref{li:mb}.  Set $f = f_{a,L}$
for the sake of simplicity.

Let $k_{1} = \frac{1}{2} \min \{ \Phi''(\psup{c}{1}), \Phi''(\psup{c}{2}) \}$.
There exists $\de_{2} = \de_{2} (\Phi)$ such that $|\psup{c}{2} - \psup{c}{1}| >
2 \de_{2}$ and $|\Phi''| > k_{1}$ on $C_{\de_{2}}$.  Notice that $|f''| > k_{1}
L$ on $C_{\de_{2}}$.  At this point we introduce the auxiliary constant $K$.
This constant will be used to bound the derivative of $f$ from below away from
the critical set.  Lower bounds on $K$ will be given as the proof develops.  We
choose $K$ before we choose $L$.  Let $\si = 2 \subsup{k}{1}{-1} L^{-1} K^{3}$
and assume $\frac{\si}{2} < \de_{2}$.  For $x \in C_{\de_{2}} \setminus
C_{\frac{1}{2} \si}$ we have $|f'(x)| \geqs K^{3}$.  Choose $L$ sufficiently
large so that $|f'| \geqs K^{3}$ outside $C_{\frac{1}{2} \si}$.  Summarizing,
the map $f$ has the following properties.

\begin{list}{\bfseries (\theprealm)}
{
\usecounter{prealm}
\setlength{\topsep}{1.5ex plus 0.2ex minus 0.2ex}
\setlength{\labelwidth}{1.2cm}
\setlength{\leftmargin}{1.5cm}
\setlength{\labelsep}{0.3cm}
\setlength{\rightmargin}{0.5cm}
\setlength{\parsep}{0.5ex plus 0.2ex minus 0.1ex}
\setlength{\itemsep}{0ex plus 0.2ex}
}
\item
$|f''| > k_{1} L$ on $C_{\de_{2}}$\label{li:p1}
\item
$|f'| \geqs K^{3}$ outside $C_{\frac{1}{2} \si}$\label{li:p2}
\end{list}

The following recovery lemma asserts that if an orbit visits a small
neighborhood of a critical point, then the derivative along this orbit regains a
definite amount of exponential growth as this orbit tracks the orbit of the
critical point for a period of time.  Set $K_{2} = \| \Phi \|_{C^{2}}$.  Let $V
= \{ x \in S^{1} : |f'(x)| \leqs K \}$ and note that $V \subset C_{\frac{1}{2}
\si}$.  Together with~\pref{li:p1} and~\pref{li:p2}, Lemma~\ref{l:recove}
implies that if $f$ satisfies Definition~\ref{d:mis}\pref{li:mb}, then $f \in
\mscr{M}$.

\begin{lemma}[Recovery estimate]\label{l:recove}
Let $c \in C$ be such that $f^{n} (c) \notin C_{\si}$ for all $n \in \NN$.  For
$x \in V$, let $n(x)$ be the smallest value of $n$ such that $|f^{n} (x) - f^{n}
(c)| > \frac{1}{4 K_{2}} K^{3} L^{-1}$.  We have $n(x) > 1$ and
$|(f^{n(x)})'(x)| \geqs k_{3} K^{n(x)}$ for some $k_{3} = k_{3} (k_{1}, K_{2})$.
\end{lemma}

\noindent
The proof of Lemma~\ref{l:recove} uses the following distortion estimate.

\begin{sublemma}[Local distortion estimate]\label{sl:lde}
Let $x$, $y \in S^{1}$.  For $i \in \ZZ^{+}$, let $\om_{i}$ denote the segment
between $f^{i} (x)$ and $f^{i} (y)$.  If $n \in \ZZ^{+}$ is such that $|\om_{i}|
\leqs \frac{1}{4 K_{2}} K^{3} L^{-1}$ and $d(\om_{i}, C) \geqs \frac{1}{2} \si$
for all $0 \leqs i < n$, then $\frac{(f^{n})'(x)}{(f^{n})'(y)} \leqs 2$.
\end{sublemma}

\begin{proof}[Proof of Sublemma~\ref{sl:lde}]
We have
\begin{align*}
\log \left( \frac{(f^{n})'(x)}{(f^{n})'(y)} \right) &= \sum_{i=0}^{n-1} \log
\left( \frac{f'(f^{i} (x))}{f'(f^{i} (y))} \right)\\
&\leqs \sum_{i=0}^{n-1} \frac{|f'(f^{i} (x)) - f'(f^{i} (y))|}{|f'(f^{i}
(y))|}\\
&\leqs \sum_{i=0}^{n-1} \frac{L K_{2} |f^{i} (x) - f^{i} (y)|}{K^{3}}\\
&\leqs \frac{L K_{2}}{K^{3}} \left( \sum_{i=0}^{n-1} \frac{1}{K^{3i}} \right)
|f^{n-1} (x) - f^{n-1} (y)| < \log (2)
\end{align*}
provided $K$ is sufficiently large.
\end{proof}

\begin{proof}[Proof of Lemma~\ref{l:recove}]

We first show that $n(x) > 1$.  Since $x \in V$, we have
\begin{equation*}
K \geqs |f'(x)| = |f''(\ga_{1})| \cdot |x-c|
\end{equation*}
and therefore
\begin{equation*}
|f(x) - f(c)| = \frac{1}{2} |f''(\ga_{2})| \cdot |x-c|^{2} \leqs
\frac{|f''(\ga_{2})|}{2 |f''(\ga_{1})|^{2}} K^{2} \leqs \frac{\| \Phi''
\|_{C^{0}}}{2 k_{1} L} K^{2}.
\end{equation*}
We may assume the final quantity is less than $\frac{1}{4 K_{2}} K^{3} L^{-1}$.
If $n(x) = 2$, then Sublemma~\ref{sl:lde} implies
\begin{equation*}
\frac{1}{4 K_{2}} K^{3} L^{-1} < |f^{2} (x) - f^{2} (c)| = |(f^{2})'(\ga_{3})|
\cdot |x-c| \leqs 2 |(f^{2})' (x)| \cdot |x-c|.
\end{equation*}
This inequality coupled with the estimate $|x-c| \leqs \frac{K}{L \| \Phi''
\|_{C^{0}}}$ implies
\begin{equation*}
|(f^{2})' (x)| > \frac{\| \Phi'' \|_{C^{0}}}{8 K_{2}} K^{2}.
\end{equation*}
Now assume $n = n(x) \geqs 3$.  Applying Sublemma~\ref{sl:lde} to estimate
$|f^{n-1} (x) - f^{n-1} (c)|$ and $|f^{n} (x) - f^{n} (c)|$, we have
\begin{gather} 
\frac{1}{2} |f''(\ga_{2})| \cdot |x-c|^{2} \cdot \frac{1}{2}
|(f^{n-2})'(f(c))| \leqs \frac{1}{4 K_{2}} K^{3} L^{-1},\label{e:re1}\\
\frac{1}{2} |f''(\ga_{2})| \cdot |x-c|^{2} \cdot 2 |(f^{n-1})'(f(c))| >
\frac{1}{4 K_{2}} K^{3} L^{-1}.\label{e:re2} 
\end{gather}
The recovery estimate follows from the lower bound
\begin{equation*}
|(f^{n})'(x)| \geqs \frac{1}{2} |f''(\ga_{1})| \cdot |x-c| \cdot
|(f^{n-1})'(f(c))|.
\end{equation*}
Replacing $|(f^{n-1})'(f(c))|$ with the lower bound provided by~\pref{e:re2} and
then replacing $|x-c|^{-1}$ with the lower bound provided by~\pref{e:re1} yields
\begin{equation*}
|(f^{n})'(x)| \geqs \frac{k_{1}}{8 K_{2}} K^{\tfrac{3}{2} n - \tfrac{3}{2}} \geqs
\frac{k_{1}}{8 K_{2}} K^{n}.
\end{equation*}
\end{proof}
We have shown that for $L$ sufficiently large, if $f$ satisfies
Definition~\ref{d:mis}\pref{li:mb}, then $f \in \mscr{M}$ with $V = \{ x \in
S^{1} : |f'(x)| \leqs K \}$.  We now find $a^{*} \in [0, 2 \pi)$ and $L^{*} \in
[L, L + 2 \pi / |\Phi (\psup{c}{2}) - \Phi (\psup{c}{1})|)$ such that $f_{a^{*},
L^{*}}$ satisfies Definition~\ref{d:mis}\pref{li:mb} by applying
Proposition~\ref{p:m}.  Additional lower bounds on $K$ will be given as the need
arises.

Hypotheses~\pref{li:h1}-\pref{li:h3} are verified as follows.  Let $z \in S^{1}$
be such that $d(z,C) \geqs d(y,C)$ for all $y \in S^{1}$.  We have
$d(z,C_{\frac{1}{2} \si}) \geqs \frac{\pi}{2} - \frac{1}{2} \si$.  There exists
$\tilde{a} \in [0, 2 \pi)$ and $\wt{L} \in [L, L + 2 \pi / |\Phi (\psup{c}{2}) -
\Phi (\psup{c}{1})|)$ such that $\subpsup{\ga}{1}{1} (\tilde{a}, \wt{L}) =
\subpsup{\ga}{1}{2} (\tilde{a}, \wt{L}) = z$.  This is so because
\begin{equation*}
|\parop{L} \subpsup{\ga}{1}{2} (a,L) - \parop{L} \subpsup{\ga}{1}{1} (a,L)| =
|\Phi (\psup{c}{2}) - \Phi (\psup{c}{1})|.
\end{equation*}
Referring to the setting of Subsections~\ref{ss:gr} and~\ref{ss:vh1}, we have
$J(a,L) = \Phi (\psup{c}{2}) - \Phi (\psup{c}{1})$ and we therefore set $k_{0} =
|\Phi (\psup{c}{2}) - \Phi (\psup{c}{1})|$.  Let $\De_{1}$ be a parameter
interval in $a$-space of length $\frac{\la_{M}}{\subsup{k}{0}{2}} K^{-3}$
centered at $\tilde{a}$ and let $\De_{2}$ be a parameter interval in $L$-space
of the same length centered at $\wt{L}$.  We assume $K$ is sufficiently large so
that $\De_{2} \subset [L, L + 2 \pi / |\Phi (\psup{c}{2}) - \Phi
(\psup{c}{1})|)$.  The image $\Ga_{1} (\De_{1} \times \De_{2})$ contains a box
such that the length of each of the sides is equal to $K^{-3}$.  Let $I_{1}$ and
$I_{2}$ be the vertical and horizontal projections of this box onto $I$,
respectively.  Since
\begin{equation*}
|\subpsup{\ga}{1}{i} (a,L) - z| \leqs \max \{ 1, |\Phi(\psup{c}{1})|, |\Phi
(\psup{c}{2})| \} \cdot \frac{\la_{M}}{2 \subsup{k}{0}{2}} K^{-3} <
\frac{\pi}{2} - \frac{1}{2} \si
\end{equation*}
for $i \in \{ 1, 2 \}$ and for all $(a,L) \in \De_{1} \times \De_{2}$ provided
$K$ is sufficiently large, we have $I_{1} \subset S^{1} \setminus C_{\frac{1}{2}
\si}$ and $I_{2} \subset S^{1} \setminus C_{\frac{1}{2} \si}$.  

By construction,~\pref{li:h1} is satisfied with $n_{0} = 1$ and the intervals
$I_{1}$ and $I_{2}$ satisfy~\pref{li:h3}.  Setting $\ve_{1} = 1$, \pref{li:h2}
is satisfied because $|f_{a,L}'| \geqs K^{3}$ on $S^{1} \setminus C_{\frac{1}{2}
\si}$ for all $(a,L) \in \De_{1} \times \De_{2}$.  If $K$ is large enough so
that~\pref{e:parset} holds, then the application of Proposition~\ref{p:m} with
$\de_{1} = \frac{1}{2} \si$ produces a Misiurewicz pair $(a^{*}, L^{*}) \in
\De_{1} \times \De_{2}$.
\end{proof}

\section{Theory of rank one attractors}\label{s:r1t}

Let $D$ denote the closed unit disk in $\RR^{n-1}$ and let $M = S^{1} \times D$.
We consider a family of maps $T_{\mbf{a}, b} : M \to M$, where $\mbf{a} =
(a_{1}, \ldots, a_{k}) \subset \Om$ is a vector of parameters and $b \in B_{0}$
is a scalar parameter.  Here $\Om = \Om_{1} \times \cdots \times \Om_{k} \subset
\RR^{k}$ is a product of intervals and $B_{0} \subset \RR \setminus \{ 0 \}$ is
a subset of $\RR$ with an accumulation point at $0$.  Points in $M$ are denoted
by $(x,y)$ with $x \in S^{1}$ and $y \in D$.  Rank one theory postulates the
following.

\begin{list}{\bfseries (\thegrealm)}
{
\usecounter{grealm}
\setlength{\topsep}{1.5ex plus 0.2ex minus 0.2ex}
\setlength{\labelwidth}{0.9cm}
\setlength{\leftmargin}{1.1cm}
\setlength{\labelsep}{0.2cm}
\setlength{\rightmargin}{0.0cm}
\setlength{\parsep}{0.5ex plus 0.2ex minus 0.1ex}
\setlength{\itemsep}{0ex plus 0.2ex}
}
\item
{\bfseries Regularity conditions.}\label{li:g1}
\begin{enumerate}[\bfseries (a)]
\item
For each $b \in B_{0}$, the function $(x,y,\mbf{a}) \mapsto T_{\mbf{a},b} (x,y)$
is $C^{3}$.\label{li:g1a}
\item
Each map $T_{\mbf{a},b}$ is am embedding of $M$ into itself.\label{li:g1b}
\item
There exists $K_{D} > 0$ independent of $\mbf{a}$ and $b$ such that for all
$\mbf{a} \in \Om$, $b \in B_{0}$, and $z$, $z' \in M$, we have\label{li:g1c}
\begin{equation*}
\frac{|\det DT_{\mbf{a},b} (z)|}{|\det DT_{\mbf{a},b} (z')|} \leqs K_{D}.
\end{equation*}
\end{enumerate}
\item
{\bfseries Existence of a singular limit.}\label{li:g2} For $\mbf{a} \in \Om$,
there exists a map $T_{\mbf{a},0} : M \to S^{1} \times \{ 0 \}$ such that the
following holds.  We select a special index $j \in \{ 1, \ldots, k \}$.  For
every fixed set $\{ a_{i} \in \Om_{i} : i \neq j \}$, the maps $(x, y, a_{j})
\mapsto T_{\mbf{a}, b} (x, y)$ converge in the $C^{3}$ topology to $(x, y,
a_{j}) \mapsto T_{\mbf{a}, 0} (x, y)$.  Identifying $S^{1} \times \{ 0 \}$ with
$S^{1}$, we refer to $T_{\mbf{a},0}$ and the restriction $f_{\mbf{a}} : S^{1}
\to S^{1}$ defined by $f_{\mbf{a}} (x) = T_{\mbf{a},0} (x,0)$ as the {\bfseries
\itshape singular limit} of $T_{\mbf{a},b}$.
\item
{\bfseries Existence of a sufficiently expanding map within the singular
limit.}\label{li:g3} 
There exists $\mbf{a}^{*} = (\subsup{a}{1}{*}, \ldots, \subsup{a}{k}{*}) \in
\Om$ such that $f_{\mbf{a}^{*}} \in \mscr{M}$.
\item
{\bfseries Parameter transversality.}\label{li:g4}
Let $C_{\mbf{a}^{*}}$ denote the critical set of $f_{\mbf{a}^{*}}$.  Define
$\tilde{\mbf{a}}_{j} = (\subsup{a}{1}{*}, \ldots, \subsup{a}{j-1}{*}, a_{j},
\subsup{a}{j+1}{*}, \ldots, \subsup{a}{k}{*})$.  We say that the family $\{
f_{\mbf{a}} \}$ satisfies the {\bfseries \itshape parameter transversality}
condition with respect to parameter $a_{j}$ if the following holds.  For each $x
\in C_{\mbf{a}^{*}}$, let $p = f(x)$ and let $x(\tilde{\mbf{a}}_{j})$ and
$p(\tilde{\mbf{a}}_{j})$ denote the continuations of $x$ and $p$, respectively,
as the parameter $a_{j}$ varies around $\subsup{a}{j}{*}$.  The point
$p(\tilde{\mbf{a}}_{j})$ is the unique point such that $p(\tilde{\mbf{a}}_{j})$
and $p$ have identical itineraries under $f_{\tilde{\mbf{a}}_{j}}$ and
$f_{\mbf{a}^{*}}$, respectively.  We have
\begin{equation*}
\left. \frac{d}{d a_{j}} f_{\tilde{\mbf{a}}_{j}} (x(\tilde{\mbf{a}}_{j}))
\right|_{a_{j} = \subsup{a}{j}{*}} \neq \left. \frac{d}{d a_{j}}
p(\tilde{\mbf{a}}_{j}) \right|_{a_{j} = \subsup{a}{j}{*}}.
\end{equation*}
\item
{\bfseries Nondegeneracy at \textquoteleft turns'.}\label{li:g5} 
For each $x \in C_{\mbf{a}^{*}}$, there exists $1 \leqs \ell \leqs n-1$ such
that
\begin{equation*}
\parfracop{y_{\ell}} T_{\mbf{a}^{*}, 0} (x,0) \neq 0.
\end{equation*}
\item
{\bfseries Conditions for mixing.}\label{li:g6}
\begin{enumerate}[\bfseries (a)]
\item
We have $e^{\tfrac{1}{3} \lambda_{0}} > 2$, where $\lambda_{0}$ is defined
within Definition~\ref{d:mis}.\label{li:g6a}
\item
\label{li:g6b}
Let $J_{1}, \ldots, J_{r}$ be the intervals of monotonicity of
$f_{\mbf{a}^{*}}$.  Let $Q = (q_{ij})$ be the matrix defined by
\begin{equation*}
q_{ij} =
\begin{cases}
1, &\text{if } f_{\mbf{a}^{*}} (J_{i}) \supset J_{j},\\
0, &\text{otherwise}.
\end{cases}
\end{equation*}
There exists $N > 0$ such that $Q^{N} > 0$.
\end{enumerate}
\end{list}

\noindent
The following lemma often facilitates the verification of~\bpref{li:g4}.

\begin{lemma}[\cite{TpTcYls1992, TpTcYls1994}]\label{l:ptreform}
Let $f = f_{\mbf{a}^{*}}$.  Suppose that for all $x \in C_{\mbf{a}^{*}}$, we
have
\begin{equation*}
\sum_{k = 0}^{\infty} \frac{1}{|(f^{k})' (f(x))|} < \infty.
\end{equation*}
Then for each $x \in C_{\mbf{a}^{*}}$,
\begin{equation}\label{e:ptreform}
\sum_{k = 0}^{\infty} \frac{[(\parop{a_{j}} f_{\tilde{\mbf{a}}_{j}}) (f^{k}
(x))]_{a_{j} = a_{j}^{*}}}{(f^{k})' (f(x))} = \left[ \frac{d}{d a_{j}}
f_{\tilde{\mbf{a}}_{j}} (x(\tilde{\mbf{a}}_{j})) - \frac{d}{d a_{j}}
p(\tilde{\mbf{a}}_{j}) \right]_{a_{j} = a_{j}^{*}}.
\end{equation}
\end{lemma}

Rank one theory states that given a family $\{ T_{\mbf{a},b} \}$
satisfying~\textbf{\pref{li:g1}-\pref{li:g5}}, a measure-theoretically
significant subset of this family consists of maps admitting attractors with
strong chaotic and stochastic properties.  We formulate the precise results and
we then describe the properties that the attractors possess.

\begin{theorem}[\cite{WqYls2001, WqYls2005}]\label{t:wyt1}
Suppose the family $\{ T_{\mbf{a},b} \}$
satisfies~\textbf{\pref{li:g1}-\pref{li:g3}} and~\bpref{li:g5}.  For all $1
\leqs j \leqs k$ such that the parameter $a_{j}$ satisfies~\bpref{li:g4} and for
all sufficiently small $b \in B_{0}$, there exists a subset $A_{j} \subset
\Om_{j}$ of positive Lebesgue measure such that for $a_{j} \in A_{j}$,
$T_{\tilde{\mbf{a}}_{j}, b}$ admits a strange attractor $\La$ with
properties~\bpref{li:sa1},~\bpref{li:sa2}, and~\bpref{li:sa3}.
\end{theorem}

\begin{theorem}[\cite{WqYls2001, WqYls2002, WqYls2005}]\label{t:wyt2}
In the sense of Theorem~\ref{t:wyt1},
\begin{equation*}
\text{\textbf{\pref{li:g1}-\pref{li:g6}}} \Longrightarrow
\text{\textbf{\pref{li:sa1}-\pref{li:sa4}}}.
\end{equation*}
\end{theorem}

\begin{list}{\bfseries (\thesarealm)}
{
\usecounter{sarealm}
\setlength{\topsep}{1.5ex plus 0.2ex minus 0.2ex}
\setlength{\labelwidth}{1.15cm}
\setlength{\leftmargin}{1.35cm}
\setlength{\labelsep}{0.2cm}
\setlength{\rightmargin}{0.0cm}
\setlength{\parsep}{0.5ex plus 0.2ex minus 0.1ex}
\setlength{\itemsep}{0ex plus 0.2ex}
}
\item
{\bfseries Positive Lyapunov exponent.}\label{li:sa1}
Let $U$ denote the basin of attraction of the attractor $\La$.  For almost every
$(x,y) \in U$ with respect to Lebesgue measure, the orbit of $(x,y)$ has a
positive Lyapunov exponent.  That is,
\begin{equation*}
\lim_{n \to \infty} \frac{1}{n} \log \| D T^{n} (x,y) \| > 0.
\end{equation*}
\item
{\bfseries Existence of SRB measures and basin property.}\label{li:sa2}
\begin{enumerate}[\bfseries (a)]
\item
The map $T$ admits at least one and at most finitely many ergodic SRB measures
all of which have no zero Lyapunov exponents.  Let $\nu_{1}, \cdots, \nu_{r}$
denote these measures.\label{li:sa2a}
\item
For Lebesgue-a.e. $(x,y) \in U$, there exists $j(x) \in \{ 1, \ldots, r \}$ such
that for every continuous function $\vp : U \to \RR$,\label{li:sa2b}
\begin{equation*}
\frac{1}{n} \sum_{i=0}^{n-1} \vp (T^{i} (x,y)) \to \int \vp \, d \nu_{j(x)}.
\end{equation*}
\end{enumerate}
\item
{\bfseries Statistical properties of dynamical observations.}\label{li:sa3}
\begin{enumerate}[\bfseries (a)]
\item
For every ergodic SRB measure $\nu$ and every H\"{o}lder continuous function
$\vp : \La \to \RR$, the sequence $\{ \vp \circ T^{i} : i \in \ZZ^{+} \}$ obeys
a central limit theorem.  That is, if $\int \vp \, d \nu = 0$, then the
sequence\label{li:sa3a}
\begin{equation*}
\frac{1}{\sqrt{n}} \sum_{i=0}^{n-1} \vp \circ T^{i}
\end{equation*}
converges in distribution to the normal distribution.  The variance of the
limiting normal distribution is strictly positive unless $\vp \circ T = \psi
\circ T - \psi$ for some $\psi$.
\item
\label{li:sa3b}
Suppose that for some $N \geqs 1$, $T^{N}$ has an SRB measure $\nu$ that is
mixing.  Then given a H\"{o}lder exponent $\eta$, there exists $\tau = \tau
(\eta) < 1$ such that for all H\"{o}lder $\vp$, $\psi : \La \to \RR$ with
H\"{o}lder exponent $\eta$, there exists $L = L(\vp, \psi)$ such that for all $n
\in \NN$,
\begin{equation*}
\left| \int (\vp \circ T^{nN}) \psi \, d \nu - \int \vp \, d \nu \int \psi \, d
\nu \right| \leqs L(\vp, \psi) \tau^{n}.
\end{equation*}
\end{enumerate}
\item
{\bfseries Uniqueness of SRB measures and ergodic properties.}\label{li:sa4}
\begin{enumerate}[\bfseries (a)]
\item
The map $T$ admits a unique (and therefore ergodic) SRB measure $\nu$,
and\label{li:sa4a}
\item
the dynamical system $(T, \nu)$ is mixing, or, equivalently, isomorphic to a
Bernoulli shift.\label{li:sa4b}
\end{enumerate}
\end{list}

\section{Proof of Theorem~\ref{t:sink}}\label{s:pst}

\subsection{Degenerate Hopf bifurcation: the reduced equations}\label{ss:dhbrs}

We study the two-dimensional system
\begin{equation}\label{e:dhbrs}
\left\{
\begin{aligned}
\dot{r} &= - \mu r\\
\dot{\thet} &= \om + \ga_{\mu} \mu + \be_{\mu} r^{2}
\end{aligned}
\right.
\end{equation}
System~\pref{e:dhbrs} is obtained from~\pref{e:dhbnf} by setting $g_{\mu} =
h_{\mu} = 0$.  Let $F_{t}$ denote the flow of~\pref{e:dhbrs}.  For $\mu$
sufficiently large, $F_{\ta (\mu)} \circ \ka$ maps $\mscr{A}$ into $\mscr{A}$.
The study of this annulus map is the central goal of this subsection.

We introduce a new coordinate system in order to standardize the position and
size of $\mscr{A}$.  Let $r = \mu^{\rh_{1}} z$.  Written in terms of $z$ and
$\thet$, system~\pref{e:dhbrs} becomes
\begin{equation}\label{e:dhbrsz}
\left\{
\begin{aligned}
\dot{z} &= - \mu z\\
\dot{\thet} &= \om + \ga_{\mu} \mu + \mu^{2 \rh_{1}} \be_{\mu} z^{2}
\end{aligned}
\right.
\end{equation}
Let $G_{t}$ denote the flow associated with~\pref{e:dhbrsz}.  The kick map $\ka$
is now given in rectangular coordinates by
\begin{equation*}
\ka
\begin{pmatrix}
z \cos (\thet)\\
z \sin (\thet)
\end{pmatrix}
=
\begin{pmatrix}
z \cos (\thet)\\
z \sin (\thet) + L \mu^{\rh_{2} - \rh_{1}}
\end{pmatrix}
\end{equation*}
We have $\mscr{A} = \{ (z, \thet) : K_{4}^{-1} \leqs z \leqs K_{4} \}$.  The
relaxation time $\ta (\mu)$ is given by
\begin{equation}\label{e:relaxz}
\tilde{z} e^{- \mu \ta (\mu)} = 1,
\end{equation}
where $\tilde{z} = L \mu^{\rh_{2} - \rh_{1}} - K_{4}$.  Let $\Psi_{\mu} = G_{\ta
(\mu)} \circ \ka$.  For $\mu$ sufficiently large, $\Psi_{\mu}$ maps $\mscr{A}$
into $\mscr{A}$.  We now derive $\Psi_{\mu} : \mscr{A} \to \mscr{A}$ explicitly.

Let $(z_{0}, \thet_{0}) \in \mscr{A}$.  Writing $\ka (z_{0}, \thet_{0}) =
(z_{1}, \thet_{1})$, we have
\begin{gather}
\subsup{z}{1}{2} = \subsup{z}{0}{2} + 2L \mu^{\rh_{2} - \rh_{1}} z_{0} \sin
(\thet_{0}) + L^{2} \mu^{2 (\rh_{2} - \rh_{1})},\label{e:z1sq}\\
\thet_{1} = \frac{\pi}{2} - \tan^{-1} \left( \frac{z_{0} \cos (\thet_{0})}{z_{0}
\sin (\thet_{0}) + L \mu^{\rh_{2} - \rh_{1}}} \right).\notag
\end{gather}
Integrating~\pref{e:dhbrsz} and writing $G_{t} \circ \ka = (z(t), \thet (t))$,
we have
\begin{align*}
z(t) &= z_{1} e^{- \mu t},\\ 
\thet (t) &= \thet_{1} + t (\om + \ga_{\mu} \mu) + \frac{\be_{\mu}}{2} \mu^{2
\rh_{1} - 1} \subsup{z}{1}{2} (1 - e^{-2 \mu t}).
\end{align*}
Evaluating $\thet (\ta (\mu))$ using~\pref{e:relaxz}, we have
\begin{equation}\label{e:thettmu1}
\thet (\ta (\mu)) = \thet_{1} + (\om + \ga_{\mu} \mu) \ta (\mu) +
\frac{\be_{\mu}}{2} \mu^{2 \rh_{1} - 1} \left( \subsup{z}{1}{2} -
\frac{\subsup{z}{1}{2}}{\tilde{z}^{2}} \right).
\end{equation}
Replacing the first occurrence of $\subsup{z}{1}{2}$ in~\pref{e:thettmu1} with
the right side of~\pref{e:z1sq}, we obtain
\begin{equation*}
\thet (\ta (\mu)) = \thet_{1} + \xi(\mu) + \frac{\be_{\mu}}{2} \left( \mu^{2
\rh_{1} - 1} \subsup{z}{0}{2} + 2L z_{0} \sin (\thet_{0}) \mu^{\rh_{1} + \rh_{2}
- 1} - \mu^{2 \rh_{1} - 1} \frac{\subsup{z}{1}{2}}{\tilde{z}^{2}} \right),
\end{equation*}
where 
\begin{equation*} 
\xi (\mu) = (\om + \ga_{\mu} \mu) \ta (\mu) + \frac{\be_{\mu}}{2} L^{2} \mu^{2
\rh_{2} - 1}.
\end{equation*}
The second component of $\Psi_{\mu}$ is given by
\begin{equation*}
z (\ta (\mu)) = \frac{z_{1}}{\tilde{z}}.
\end{equation*}

We wish to show that the family $\{ \Psi_{\mu} \}$ converges to a singular limit
as $\mu \to 0$.  This cannot be accomplished directly because $\xi (\mu)$
diverges as $\mu \to 0$, preventing the convergence of $\thet (\ta (\mu))$.  We
overcome this difficulty by taking advantage of the fact that $\thet (t)$ is
computed modulo $2 \pi$.  Assume that $\om > 0$.  For $\mu$ sufficiently small,
$\xi (\mu)$ is monotone.  In addition, $\xi (\mu) \to \infty$ as $\mu \to 0$.
Let $(\mu_{n})$ be a sequence such that $\mu_{n} \to 0$ monotonically, $\xi$ is
monotone on $(0, \mu_{1}]$, and $\xi (\mu_{n}) \in 2 \pi \ZZ$ for all $n \in
\NN$.  We introduce the parameter $a \in [0, 2 \pi)$ and write $\Psi_{\mu}$ in
terms of $a$.  For $n \in \NN$ and $a \in [0, 2 \pi)$, let $\mu (a, n) =
\xi^{-1} (\xi (\mu_{n}) + a)$.  When referring to $\mu (a,n)$, we will
henceforth suppress the dependence on $n$ and simply write $\mu (a)$.  The
problematic term $\xi (\mu)$ becomes $\xi (\mu (a)) = a$.  Writing $\Psi_{\mu
(a)} = T_{a, L, \mu_{n}}$, we have
\begin{gather}
T_{a, L, \mu_{n}}^{\langle 1 \rangle} (z_{0}, \thet_{0}) =
\frac{z_{1}}{\tilde{z}},\notag\\
\begin{aligned}
T_{a, L, \mu_{n}}^{\langle 2 \rangle} (z_{0}, \thet_{0}) = \thet_{1} + a &+
\frac{\be_{\mu (a)}}{2} \bigg( \mu (a)^{2 \rh_{1} - 1} z_{0}^{2} +\\ 
&2 L \mu (a)^{\rh_{1} + \rh_{2} - 1} z_{0} \sin (\thet_{0}) - \mu (a)^{2 \rh_{1}
- 1} \frac{z_{1}^{2}}{\tilde{z}^{2}} \bigg),
\end{aligned}
\notag
\end{gather}
where $T^{\langle 1 \rangle}$ and $T^{\langle 2 \rangle}$ are the components of
$T$.  Let $\rh_{1}$ and $\rh_{2}$ satisfy $\frac{1}{2} < \rh_{1} < 1$ and
$\rh_{1} + \rh_{2} = 1$.  Then as $n \to \infty$, $T_{a, L, \mu_{n}}$ converges
in the $C^{0}$ topology to the map $T_{a, L, 0}$ defined by
\begin{align}
T_{a, L, 0}^{\langle 1 \rangle} &= 1\notag\\
T_{a, L, 0}^{\langle 2 \rangle} &= \frac{\pi}{2} + \be_{0} L z_{0} \sin
(\thet_{0}) + a.\notag
\end{align}
The following lemma asserts that the convergence is strong enough for the
application of rank one theory. 

\begin{lemma}\label{l:c3r}
Fix $L > 0$.  The maps $(z_{0}, \thet_{0}, a) \mapsto T_{a, L, \mu_{n}} (z_{0},
\thet_{0})$ converge in the $C^{3}$ topology to the map $(z_{0}, \thet_{0}, a)
\mapsto T_{a, L, 0} (z_{0}, \thet_{0})$ as $n \to \infty$ on the domain
$\mscr{A} \times [0, 2 \pi)$.
\end{lemma}

\begin{proof}[Proof of Lemma~\ref{l:c3r}]
Holding $a$ fixed, the derivatives of $\thet_{1}$, $\frac{z_{1}}{\tilde{z}}$,
and $\frac{z_{1}^{2}}{\tilde{z}^{2}}$ of orders $1$, $2$, and $3$ with respect
to $z_{0}$ and $\thet_{0}$ are $\mcal{O} (\mu^{\rh_{1} - \rh_{2}})$.  When
differentiating with respect to $a$, use the fact that for $i=1,2,3$,
\begin{equation*}
\parop{a}^{(i)} \mu (a) = \mcal{O} \left( \frac{\mu_{n}^{i+1}}{\big( \log
(\mu_{n}^{-1}) \big)^{i}} \right).
\end{equation*}
\end{proof}

\noindent
We finish this subsection with a distortion estimate.

\begin{lemma}[Distortion estimate]\label{l:rsde}
Let $0 < L_{2} < L_{3}$.  There exists $K_{D} > 0$ such that for all $n \in
\NN$, $a \in [0, 2 \pi)$, $L \in [L_{2}, L_{3}]$, and $(z_{0}, \thet_{0})$,
$(z_{0}', \thet_{0}') \in \mscr{A}$, we have
\begin{equation*}
\frac{| \det DT_{a, L, \mu_{n}} (z_{0}, \thet_{0}) |}{| \det DT_{a, L, \mu_{n}}
(z_{0}', \thet_{0}') |} \leqs K_{D}.
\end{equation*}
\end{lemma}

\begin{proof}[Proof of Lemma~\ref{l:rsde}]
Recall that $T_{a, L, \mu_{n}} = G_{\ta (\mu (a))} \circ \ka$.  We bound the
distortion by analyzing $G$ and $\ka$ independently.  Let $(z_{0}, \thet_{0})$,
$(z_{0}', \thet_{0}') \in \mscr{A}$.  Writing $\mu = \mu (a)$, $\det D \ka
(z_{0}, \thet_{0})$ is given by
\begin{equation*}
\frac{z_{0}^{3} + \mu^{\rh_{2} - \rh_{1}} \cdot L z_{0} (1 + z_{0})\sin 
(\thet_{0}) + \mu^{2 (\rh_{2} - \rh_{1})} \cdot L^{2} (1 + (z_{0} - 1) \cos^{2} 
(\thet_{0}))}{z_{1} \big( (z_{0} \sin (\thet_{0}) + L \mu^{\rh_{2} - 
\rh_{1}})^{2} + z_{0}^{2} \cos^{2} (\thet_{0}) \big)}.
\end{equation*}
This explicit formula implies the estimate
\begin{equation*}
\frac{\det D \ka (z_{0}, \thet_{0})}{\det D \ka (z_{0}', \thet_{0}')} = \mcal{O}
(1).
\end{equation*}
Now set $G = G_{\ta (\mu (a))}$.  For any point in $\ka (\mscr{A})$, the
determinant of the derivative of $G$ is precisely $\tilde{z}^{-1}$.  Therefore,
\begin{equation*}
\frac{\det DG (z_{1}, \thet_{1})}{\det DG (z_{1}', \thet_{1}')} = 1.
\end{equation*}
\end{proof}

\subsection{Inclusion of the higher-order terms in the normal
form}\label{ss:dhbnf}

We show that the inclusion of the higher-order terms in the differential
equations defining the flow does not affect the form of the singular limit
derived in Subsection~\ref{ss:dhbrs}.  Set $r = \mu^{\rh_{1}} \hat{z}$ and
$\thet = \hat{\thet}$.  Written in terms of $\hat{z}$ and $\hat{\thet}$, the
normal form~\pref{e:dhbnf} becomes
\begin{equation}\label{e:dhbnfz}
\left\{
\begin{aligned}
\dot{\hat{z}} &= - \mu \hat{z} + \mu^{4 \rh_{1}} \hat{z}^{5} g_{\mu}
(\mu^{\rh_{1}} \hat{z}, \hat{\thet})\\
\dot{\hat{\thet}} &= \om + \ga_{\mu} \mu + \be_{\mu} \mu^{2 \rh_{1}} \hat{z}^{2}
+ \mu^{4 \rh_{1}} \hat{z}^{4} h_{\mu} (\mu^{\rh_{1}} \hat{z}, \hat{\thet})
\end{aligned}
\right.
\end{equation}
Let $\wh{G}_{t}$ denote the flow generated by~\pref{e:dhbnfz}.  We define the
family $\{ \wh{T} \}$ on $\mscr{A}$ by first applying the kick map $\ka$ and
then allowing the $\wh{G}_{t}$-flow to return $\ka (\mscr{A})$ to $\mscr{A}$.
Set $\wh{T}_{a, L, \mu_{n}} = \wh{G}_{\ta (\mu (a))} \circ \ka$.

\begin{lemma}\label{l:c3nfk}
Fix $L > 0$.  If $\rh_{2} > \frac{1}{3}$, then the maps $(z_{0}, \thet_{0}, a)
\mapsto \wh{T}_{a, L, \mu_{n}} (z_{0}, \thet_{0})$ converge in the $C^{3}$
topology to the map $(z_{0}, \thet_{0}, a) \mapsto T_{a, L, 0} (z_{0},
\thet_{0})$ as $n \to \infty$ on the domain $\mscr{A} \times [0, 2 \pi)$.
\end{lemma}

\begin{proof}[Proof of Lemma~\ref{l:c3nfk}]
Computing the first component of $\wh{G}_{t} \circ \ka$, we have
\begin{align}
\hat{z} (t) &= z_{1} e^{-\mu t} \left( 1 + \mu^{4 \rh_{1}} \int_{0}^{t}
z_{1}^{-1} e^{\mu s} \hat{z} (s)^{5} g_{\mu} (\mu^{\rh_{1}} \hat{z} (s),
\hat{\thet} (s)) \, ds \right)\notag\\
&= z(t) + \ze (t),\notag
\end{align}
where the perturbative term $\ze (t)$ is defined by
\begin{equation}\label{e:com1pt}
\ze (t) = \mu^{4 \rh_{1}} e^{-\mu t} \int_{0}^{t} e^{\mu s} \hat{z} (s)^{5} 
g_{\mu} (\mu^{\rh_{1}} \hat{z} (s), \hat{\thet} (s)) \, ds.
\end{equation}
Computing the second component of $\wh{G}_{t} \circ \ka$, we have $\hat{\thet}
(t) = \thet (t) + \tilde{\thet} (t)$, where
\begin{equation}\label{e:com2pt}
\begin{aligned}
\tilde{\thet} (t) = &\be_{\mu} \mu^{2 \rh_{1}} \int_{0}^{t} \bigg[ 2 \mu^{4 
\rh_{1}} z_{1} e^{-2 \mu v} \int_{0}^{v} e^{\mu s} \hat{z} (s)^{5} g_{\mu} 
(\mu^{\rh_{1}} \hat{z} (s), \hat{\thet} (s)) \, ds\\
&+ \mu^{8 \rh_{1}} e^{-2 \mu v} \left( \int_{0}^{v} e^{\mu s} \hat{z} (s)^{5}
g_{\mu} (\mu^{\rh_{1}} \hat{z} (s), \hat{\thet} (s)) \, ds \right)^{2} \bigg] \,
dv\\ 
&+ \mu^{4 \rh_{1}} \int_{0}^{t} \hat{z} (s)^{4} h_{\mu} (\mu^{\rh_{1}} \hat{z}
(s), \hat{\thet} (s)) \, ds.
\end{aligned}
\end{equation}
In order to establish $C^{0}$ convergence, it suffices to show that the
perturbative terms $\ze (\ta (\mu (a)))$ and $\tilde{\thet} (\ta (\mu (a)))$
converge to $0$ in the $C^{0}$ topology as $n \to \infty$.  Estimating the
integrals in~\pref{e:com1pt} and~\pref{e:com2pt}, we obtain
\begin{gather}
\ze (\ta (\mu)) = \mcal{O} (\mu^{5 \rh_{2} - \rh_{1} - 1} \log
(\mu^{-1})),\notag\\ 
\tilde{\thet} (\ta (\mu)) = \mcal{O} \big( \mu^{6 \rh_{2} - 2} (\log 
(\mu^{-1}))^{2} \big) + \mcal{O} \big( \mu^{10 \rh_{2} - 3} (\log
(\mu^{-1}))^{3} \big).\notag
\end{gather}
Since $\rh_{2} \in (\frac{1}{3}, \frac{1}{2})$ and $\rh_{1} \in (\frac{1}{2},
\frac{2}{3})$, we have
\begin{equation*}
\| \ze (\ta (\mu)) \|_{C^{0}} \to 0 \text{ and } \| \tilde{\thet} (\ta (\mu))
\|_{C^{0}} \to 0
\end{equation*}
as $\mu \to 0$.

We complete the proof of Lemma~\ref{l:c3nfk} by showing that
\begin{equation*}
\| D^{i} \wh{G}_{\ta (\mu (a))} \circ \ka - D^{i} G_{\ta (\mu (a))} \circ \ka
\|_{C^{0}} \to 0
\end{equation*}
for $1 \leqs i \leqs 3$.  In light of Lemma~\ref{l:c3r}, this establishes the
asserted $C^{3}$ convergence.  Since $\| D^{i} \ka \|_{C^{0}}$ is bounded for $1
\leqs i \leqs 3$, it is sufficient to show that $\| D^{i} \wh{G}_{\ta (\mu (a))}
- D^{i} G_{\ta (\mu (a))} \|_{C^{0}} \to 0$ for $1 \leqs i \leqs 3$.  We use the
following elementary Gronwall-type lemma.

\begin{lemma}[\cite{WqYls2003}]\label{l:gtype}
Let $\La \subset \RR^{N}$ be a convex open domain.  Let $W$ and $\wh{W}$ be
$C^{1}$ vector fields on $\La$.  Suppose that for $t \in [0, t_{0}]$,
$\hat{\vp}$ and $\vp$ solve the equations
\begin{equation*}
\frac{d \hat{\vp}}{dt} = \wh{W} (\hat{\vp}) \text{ and } \frac{d \vp}{dt} = W
(\vp)
\end{equation*}
with $\hat{\vp} (0) = \vp (0)$.  Then for all $t \in [0, t_{0}]$, we have
\begin{equation*}
\| \hat{\vp} (t) - \vp (t) \| \leqs \frac{A_{1}}{A_{2}} (e^{A_{2} t} - 1),
\end{equation*}
where
\begin{equation*}
A_{1} = \sup_{x \in \La} \| \wh{W} (x) - W(x) \| \text{ and } A_{2} =
\sum_{j=1}^{N} \sup_{x \in \La} \| DW^{\langle j \rangle} (x) \|.
\end{equation*}
\end{lemma}
We rescale time in~\pref{e:dhbrsz} and~\pref{e:dhbnfz} by setting $t = t' \ta
(\mu (a))$.  Let $\eta$ and $\hat{\eta}$ denote the rescaled vector fields.  We
have
\begin{gather}
\left\{
\begin{aligned}
\eta^{\langle 1 \rangle} &= \ta (\mu (a)) (- \mu z)\\
\eta^{\langle 2 \rangle} &= \ta (\mu (a)) (\om + \ga_{\mu} \mu + \mu^{2 \rh_{1}}
\be_{\mu} z^{2})
\end{aligned}
\right.
\notag\\
\left\{
\begin{aligned}
\hat{\eta}^{\langle 1 \rangle} &= \ta (\mu (a)) (- \mu \hat{z} + \mu^{4 \rh_{1}}
\hat{z}^{5} g_{\mu} (\mu^{\rh_{1}} \hat{z}, \hat{\thet}))\\
\hat{\eta}^{\langle 2 \rangle} &= \ta (\mu (a)) (\om + \ga_{\mu} \mu + \be_{\mu}
\mu^{2 \rh_{1}} \hat{z}^{2} + \mu^{4 \rh_{1}} \hat{z}^{4} h_{\mu} (\mu^{\rh_{1}}
\hat{z}, \hat{\thet}))
\end{aligned}
\right.
\notag
\end{gather}
We explicitly treat the case $i = 1$.  The cases $i = 2$ and $i = 3$ are handled
using the same technique.  Apply Lemma~\ref{l:gtype} with $\hat{\vp} = D
\wh{G}$, $\vp = DG$, $\wh{W} = D \hat{\eta}$, $W = D \eta$, and $t = 1$.  The
quantity $A_{2}$ is bounded.  Therefore, the estimate $A_{1} = \mcal{O} (\mu^{5
\rh_{2} - \rh_{1} - 1} \log (\mu^{-1}))$ implies that
\begin{equation}\label{e:vareq1}
\| D \wh{G}_{1} - D G_{1} \|_{C^{0}} = \mcal{O} (\mu^{5 \rh_{2} - \rh_{1} - 1}
\log (\mu^{-1})).
\end{equation}
\end{proof}

\subsection{Verification of~\pref{li:g1}-\pref{li:g6}}\label{ss:vg1g6}

Theorem~\ref{t:sink} follows from an application of Theorems~\ref{t:wyt1}
and~\ref{t:wyt2}.  Statements~\pref{li:sink1} and~\pref{li:sink2} of
Theorem~\ref{t:sink} require the verification of~\bpref{li:g1}-\bpref{li:g5} for
the family $\{ \wh{T}_{a, L, \mu_{n}} \}$.  Statement~\pref{li:sink3} of
Theorem~\ref{t:sink} requires the additional verification of~\bpref{li:g6}.

We proceed with the verification of statement~\pref{li:sink1} of
Theorem~\ref{t:sink}.  Properties~\bpref{li:g1}\bpref{li:g1a}
and~\bpref{li:g1}\bpref{li:g1b} follow from the general theory of ordinary
differential equations.  For~\bpref{li:g1}\bpref{li:g1c}, it suffices to show
that the distortion of $\wh{G}_{\ta (\mu (a))}$ is bounded because the
distortion of $\ka$ is bounded.  Using~\pref{e:vareq1}, we have
\begin{equation*}
D \wh{G}_{\ta (\mu (a))} (z_{1}, \thet_{1}) =
\begin{pmatrix}
\tilde{z}^{-1} + \ve_{1} & \ve_{2}\\
\frac{\be_{\mu}}{2} \mu^{2 \rh_{1} - 1} \left( 2z_{1} -
\frac{2z_{1}}{\tilde{z}^{2}} \right) + \ve_{3} & 1 + \ve_{4}
\end{pmatrix}
\end{equation*}
where $\ve_{j} = \mcal{O} (\mu^{5 \rh_{2} - \rh_{1} - 1} \log (\mu^{-1}))$ for
$1 \leqs j \leqs 4$.  Since $\rh_{2} > \frac{3}{8}$,
$\frac{\ve_{j}}{\tilde{z}^{-1}} \to 0$ as $\mu \to 0$ for $1 \leqs j \leqs 4$.
Therefore, we have 
\begin{equation*}
\det (D \wh{G}_{\ta (\mu (a))} (z_{1}, \thet_{1})) = \tilde{z}^{-1} + \mcal{O} 
(\mu^{5 \rh_{2} - \rh_{1} - 1} \log (\mu^{-1})).
\end{equation*}
This estimate implies that the distortion of $\wh{G}_{\ta (\mu (a))}$ is
bounded.

Lemma~\ref{l:c3nfk} establishes~\bpref{li:g2}.  Let $f_{a, L}$ denote the
restriction of $T_{a, L, 0}^{\langle 2 \rangle}$ to the circle $S^{1} = \{
(z_{0}, \thet_{0}) : z_{0} = 1 \}$.  We have
\begin{equation*}
f_{a, L} (\thet) = \frac{\pi}{2} + \be_{0} L \sin (\thet) + a.
\end{equation*}
Applying Theorem~\ref{t:al} with $\Phi (\thet) = \sin (\thet)$, $\psup{c}{1} =
\frac{\pi}{2}$, and $\psup{c}{2} = \frac{3 \pi}{2}$, if $L$ is sufficiently
large then there exist $L^{*} \in [L, L + \frac{\pi}{| \be_{0} |}]$ and $a^{*}
\in [0, 2 \pi)$ such that $f_{a^{*}, L^{*}} \in \mscr{M}$.  This
is~\bpref{li:g3}.  We establish parameter transversality~\bpref{li:g4} by
applying Lemma~\ref{l:ptreform}.  Write $f = f_{a^{*}, L^{*}}$ and $f_{a} =
f_{a, L^{*}}$.  We have $\parop{a} f_{a} (\cdot) = 1$ and $|(f^{k})' (f(x))|
\geqs K^{k}$.  Therefore, the absolute value of the left side
of~\pref{e:ptreform} is bounded below by $1 - \sum_{k = 1}^{\infty} K^{-k}$.
This quantity is positive if $K > 2$.  For~\bpref{li:g5}, observe that
\begin{equation*}
\parop{z_{0}} T_{a, L, 0}^{\langle 2 \rangle} (1, \psup{c}{1}) = \be_{0} L \neq
0 \text{ and } \parop{z_{0}} T_{a, L, 0}^{\langle 2 \rangle} (1, \psup{c}{2}) =
- \be_{0} L \neq 0.
\end{equation*}
This completes the verification of statement~\pref{li:sink1} of
Theorem~\ref{t:sink}.

Statement~\pref{li:sink2} of Theorem~\ref{t:sink} follows from the fact that for
all $L$ sufficiently large, $f_{a, L} \in \mscr{M}$ for a $\mcal{O}
(L^{-1})$-dense set of values of $a$.  Wang and Young~\cite{WqYls2002} prove
this result in a slightly different context.  The proof for the family $\{ f_{a,
L} \}$ is essentially the same.

Statement~\pref{li:sink3} of Theorem~\ref{t:sink} requires the verification of
the conditions for mixing~\bpref{li:g6}.  Property~\bpref{li:g6}\bpref{li:g6a}
holds provided $e^{\la_{0}} = K > 8$.  Property~\bpref{li:g6}\bpref{li:g6b} is
satisfied with $N = 1$ provided $L$ is sufficiently large.

\bibliographystyle{amsplain}
\bibliography{pkhopf} 

\end{document}